\documentclass[12pt,oneside]{amsart}

\newtheorem{thm}{Theorem}[section] 
\newtheorem{cor}[thm]{Corollary}
\newtheorem{example}[thm]{Example}
\newtheorem{corollary}[thm]{Corollary}
\newtheorem{lem}[thm]{Lemma} 
\newtheorem{lemma}[thm]{Lemma}
\newtheorem{prop}[thm]{Proposition}
\newtheorem{proposition}[thm]{Proposition}
\theoremstyle{definition}
\theoremstyle{remark}
\newtheorem{rem}[thm]{Remark} 
\theoremstyle{proof}

\numberwithin{equation}{section}

\newcommand{\norm}[1]{\left\Vert#1\right\Vert}
\newcommand{\snorm}[1]{\Vert#1\Vert}
\newcommand{\abs}[1]{\left\vert#1\right\vert}
\newcommand{\set}[1]{\left\{#1\right\}}
\newcommand{\brac}[1]{\left(#1\right)}
\newcommand{\scalar}[1]{\left \langle #1 \right \rangle}

\newcommand{\Real}{\mathbb{R}}

\newcommand{\eps}{\varepsilon}

\newcommand{\E}{\mathcal{E}}

\def \RR {\mathbb R}

\def \eps {\varepsilon}

\newcommand{\Vol}[1]{\textnormal{Vol} \left(#1 \right)}
\newcommand{\vol}[1]{\left | #1 \right|}
\newcommand{\VolSet}[1]{\textnormal{Vol} \set{#1}}
\newcommand{\volRatio}[2]{\brac{\frac{\vol{#1}}{\vol{#2}}}^{1/n}}
\newcommand{\VolRatio}[2]{\volRatio{#1}{#2}}

\newcommand{\VolRad}[1]{\textit{Vol.rad.}\brac{#1}}

\begin{document}

\title{On volume distribution in $2$-convex bodies}

\author{Bo'az Klartag  $\;\;\;\;$  Emanuel Milman}

\thanks{The first named author is a Clay Research Fellow and is also supported by NSF grant DMS-0456590. The second named author is supported in part by BSF and ISF}

\maketitle

\begin{abstract}
We consider convex sets whose modulus of convexity is uniformly
quadratic. First, we observe several interesting relations between
different positions of such ``$2$-convex'' bodies; in particular,
the isotropic position is a finite volume-ratio position for these
bodies. Second, we prove that high dimensional $2$-convex bodies
posses one-dimensional marginals that are approximately Gaussian.
Third, we improve for $1<p\leq 2$ some bounds on the isotropic
constant of quotients of subspaces of $L_p$ and $S_p^m$, the
Schatten Class space.
\end{abstract}


\section{Introduction}

The purpose of this note is to collect several interesting facts
related to the distribution of volume in high dimensional
$2$-convex bodies. Suppose that $K \subset \Real^n$ is a
centrally-symmetric (i.e. $K = -K$) convex body (i.e. a convex,
compact set with non-empty interior). Let $\| \cdot \|_K$ be the
norm on $\Real^n$ whose unit ball is $K$. The modulus of
convexity of $K$ is the function:
\begin{equation} \label{eqn_1242}
\delta_K(\eps) = \inf \left \{ 1 - \left \| \frac{x+y}{2} \right
\|_K \ ; \  \|x\|_K, \|y\|_K \leq 1, \| x - y \|_K \geq \eps
\right \},
\end{equation}
defined  for $0 < \eps \leq 2$. We say that $K$ is ``$2$-convex
with constant $\alpha$'' (see, e.g. \cite[Chapter
1.e]{LT-Book-II}), if for all $0 < \eps \leq 2$,
\begin{equation} \label{eqn_1244}
\delta_K(\eps) \geq \alpha \eps^2.
\end{equation}
Note that this should not be confused with the notions of
$p$-convexity or $q$-concavity (e.g. \cite[Chapter
1.d]{LT-Book-II}) defined for Banach lattices. Being $2$-convex
with constant $\alpha$ is a linearly invariant property.
Furthermore, as is evident from the definitions, if $K$ is
$2$-convex with constant $\alpha$, so is $K \cap E$ for any
subspace $E$. Thus sections of a convex body inherit the
$2$-convexity properties of the body. The same holds for
projections (see, e.g. Lemma \ref{lem:projections} below). A
basic example of $2$-convex bodies are unit balls of $L_p$ spaces
for $1<p\leq 2$, in which case $\alpha$ is of the order of $p-1$
(e.g. \cite[Chapter 1.e]{LT-Book-II}).
Consequently, also sections, projections, and sections of
projections of $L_p$-balls are $2$-convex bodies, with constants
that depend solely on $p$.

\medskip

It is well-known that the uniform measure on a $2$-convex body is
``well behaved'', in many senses (see, e.g. \cite{Gromov-Milman}
\cite{SchmuckyUniformlyConvexBodies} and \cite{BobkovLedoux}).
Questions on distribution of mass in high-dimensional convex sets
regained some interest in the last few years, and some partial
progress was obtained. We approach the study of mass distribution
in $2$-convex sets, in view of these developments. Arguably, the
most basic question regarding volume distribution in
high-dimensional convex sets is the Slicing Problem, or Hyperplane
Conjecture. This question asks whether for any convex body $K
\subset \Real^n$ of volume one, there exists a hyperplane $H
\subset \Real^n$ such that $\Vol{K \cap H} > c$, for some
universal constant $c > 0$. Here and henceforth, $\Vol{A}$ or
$|A|$ for short, denotes the volume of $A \subset \Real^n$ in its
affine hull. In the category of $2$-convex bodies, a positive
answer to this question was provided by Schmuckenschl\"ager
\cite{SchmuckyUniformlyConvexBodies}. We provide a more direct
approach to Schmuckenschl\"ager's result, that is based on an
argument of \cite{ABV}.

\begin{proposition} \label{prop:hyperplane}
Let $K \subset \Real^n$ be a centrally-symmetric convex body of
volume one. Suppose $K$ is $2$-convex with constant $\alpha$. Then
there exists a hyperplane $H \subset \Real^n$ such that:
\[
\Vol{K \cap H} \geq c \sqrt{\alpha},
\]
where $c > 0$ is a universal constant.
\end{proposition}

A centrally-symmetric convex $K \subset \Real^n$ of volume one is
said to be \emph{isotropic} or \emph{in isotropic position}, if
for any $\theta \in \Real^n$:
\[
\int_K \langle x, \theta \rangle^2 dx = L_K |\theta|^2,
\]
where $L_K$ is some quantity, independent of $\theta$, and
$|\cdot|$ is the Euclidean norm. In that case, the
\emph{isotropic constant of $K$} is defined as $L_K$. It is well
known (see, e.g. \cite{Milman-Pajor-LK}) that for any
centrally-symmetric convex $K \subset \Real^n$, there exists a
linear transformation such that $\tilde{K} = T(K)$ is isotropic.
Moreover, this map $T$ is unique up to orthogonal
transformations. We therefore define the isotropic constant of an
arbitrary centrally-symmetric convex body $K \subset \Real^n$, to
be $L_K = L_{\tilde{K}}$, where $\tilde{K}$ is an isotropic
linear image of $K$. An observation that goes back to Hensley
\cite{Hensley}, is that when K is isotropic, for any hyperplane H
through the origin:
\[
\frac{c_1}{L_K} \leq \Vol{K \cap H} \leq \frac{c_2}{L_K},
\]
where $c_1, c_2 > 0$ are universal constants. Based on this, the
Slicing Problem may be reformulated as follows (e.g.
\cite{Milman-Pajor-LK}): Is it true that for any dimension n and
any centrally-symmetric convex body $K \subset \Real^n$, we have
that $L_K \leq C$, where $C > 0$ is a universal constant?

\medskip

As a by-product of our methods, we improve for $1<p\leq 2$ a bound
for the isotropic constant of the unit balls of quotients of
subspaces of $L_p$, and establish the same bound for arbitrary
quotients of subspaces of $l_p$-Schatten-Class spaces of $m$ by
$m$ matrices, denoted $S_p^m$ (see Section 3 for definitions).
For a Banach Space $X$, we denote by $SQ_n(X)$ the family of all
centrally-symmetric convex bodies $K \subset \Real^n$, such that
$K$ is the unit ball of some subspace of a quotient of $X$.

\begin{proposition}  \label{prop:Junge-And-Schatten}
Let $1<p\leq 2$, let $X = L_p$ or $X = S_p^m$, and suppose that
$K \in SQ_n(X)$. Then,
\begin{equation}
L_K \leq C \sqrt{q} \label{eqn_1152}
\end{equation} where $q = p^* = p/(p-1)$ and $C > 0$ is a universal
constant.
\end{proposition}
Junge \cite{Junge-slicing-problem-for-quotients-of-Lp} has
previously proven a version of (\ref{eqn_1152}) with $q$ in place
of $\sqrt{q}$ for $X = L_p$. For $X = S_p^m$ and $1 \leq p \leq
2$, a universal bound on $L_K$ was established in
\cite{Schatten-LK-Bounded} when $K$ is the unit ball of $X$, and
in \cite{Guedon-Paouris-Schatten} when $K$ is the unit ball of
certain specific subspaces of $X$.

\medskip

In addition to the isotropic position, there are several other
important Euclidean structures that are associated with a given
convex body, such as John's position, minimal mean-width position,
$\ell$-position, (regular) $M$-position, etc. The relations
between these various positions in general are not clear. See
\cite{BKM-symmetrizations} for an equivalence of the hyperplane
conjecture to a certain putative relation between the isotropic
position and $M$-position. However, in the class of $2$-convex
bodies, the following holds:

\begin{proposition} \label{prop:finite-vr}
Let $K \subset \Real^n$ be a $2$-convex body with constant
$\alpha$ and of volume 1. If $K$ is in isotropic position then:
\[
c \sqrt{\alpha} \sqrt{n} D_n \subset K,
\]
where $D_n$ is the unit Euclidean ball in $\Real^n$ and $c > 0$ is
a universal constant.
\end{proposition}

That is, the isotropic position of a $2$-convex body is a
\emph{finite volume-ratio position}. The volume-ratio of a
centrally-symmetric convex body $K \subset \RR^n$ is defined as:
\begin{equation} \label{eq:vr-defn}
v.r.(K) = \min_{\E \subset K} \left( \frac{|K|}{|\E|}
\right)^{\frac{1}{n}},
\end{equation}
where the minimum runs over all ellipsoids that are contained in
$K$. If $v.r.(K) < C$, for some universal constant $C$, it is
customary to say that $K$ is a finite volume-ratio body. When the
minimum over all \emph{Euclidean balls} is bounded by a universal
constant, we will say that $K$ is in a finite volume-ratio
\emph{position}. Note that $c_1 < |\sqrt{n} D_n|^{1/n} < c_2$ for
some universal constants $c_1, c_2 > 0$, so Proposition
\ref{prop:finite-vr} implies that the isotropic position is a
finite volume-ratio position.

This conclusion is clearly false for general convex bodies, even
for convex bodies whose distance to the Euclidean ball is
universally bounded (see the example after Lemma \ref{lem1}
below). In Section 4 we establish further rigid relations between
various positions of $2$-convex bodies, that cannot hold for
arbitrary convex bodies. In particular, recall that $K$ is said
to be in John's maximal-volume ellipsoid position when the
minimum in (\ref{eq:vr-defn}) is attained by a Euclidean ball. We
will see the following:

\begin{prop} \label{prop:inverse-finite-vr}
Let $K \subset \Real^n$ be a $2$-convex body with constant
$\alpha$ and of volume 1. If $K$ is in John's maximal-volume
ellipsoid position, then:
\begin{equation} \label{eq:almost-iso}
\brac{\int_K |x|^2 dx}^{\frac{1}{2}} \leq \frac{C}{\alpha}
\sqrt{n},
\end{equation}
where $C>0$ is a universal constant.
\end{prop}

The latter is in a sense a converse to Proposition
\ref{prop:finite-vr}, since (\ref{eq:almost-iso}) implies that
$K$ is ``essentially" isotropic. To see this, note (e.g.
\cite{Milman-Pajor-LK}) that the isotropic position minimizes the
value of $\int_{T(K)} |x|^2 dx$, over all volume 1 affine images
$T(K)$ of $K$, and in that case we have:
\[
\inf \brac{\int_{T(K)} |x|^2 dx}^{\frac{1}{2}} = \sqrt{n} L_K.
\]
In addition to being an ``essentially" isotropic position, we show
in Section 4 that John's position is in fact an ``essentially"
minimal mean-width position and a 2-regular M-position (see
Section 4 for definitions). A complete list of other relations
between the aforementioned various positions is given at the end
of Section 4.

\medskip

An additional interesting volumetric question, is the so-called
``Central Limit Property of Convex Bodies''. Let $X$ denote a
uniformly distributed vector inside a convex set $K \subset
\Real^n$ of volume one. In its weakest form, a conjecture of
Antilla, Ball and Perissinaki \cite{ABP}, states that for some
non-zero vector $\theta \in \Real^n$, the random variable
$\langle X, \theta \rangle$ is very close to a Gaussian random
variable. That is, the total variation distance between the
random variable $\langle X, \theta \rangle$ and a corresponding
Gaussian random variable, is smaller than $\eps_n$, where
$\eps_n$ is a sequence tending to zero, that depends solely on
$n$. In this note, we verify the following (see Theorem
\ref{thm:2-convex-Gaussian} for an exact formulation):
\begin{prop} \label{prop:CLP}
The ``Central Limit Property'' holds true for arbitrary $2$-convex
bodies.
\end{prop}
In \cite{ABP}, the existence of approximately Gaussian marginals
of $2$-convex bodies was proven only under a certain, rather
weak, constraint on the diameter of $K$ in isotropic position. We
show in Example \ref{example} that there exist $2$-convex bodies
in $\Real^n$ for which this constraint is violated. In fact, we
show that there exist such bodies of volume 1 whose diameter in
isotropic position is greater than $c n$ (where $c>0$ is a
universal constant). Our idea is to put $K$ in another position,
namely L\"{o}wner's minimal diameter position, in which we show
in Proposition \ref{prop:small-diam} that the diameter is not
larger than $\frac{C}{\lambda} n^{1-\lambda}$, where $\lambda$
depends only on $\alpha$, the $2$-convexity constant of $K$ and
$C>0$ is a universal constant. We conclude Proposition
\ref{prop:CLP} by proving a version of a Theorem from \cite{ABP}
about the existence of Gaussian marginals, where the assumption
of being in isotropic position is removed (see Theorem
\ref{thm:extended-ABP}). Further developments on the existence of
Gaussian marginals of uniformly convex bodies are discussed in
\cite{EMilman-Gaussian-Marginals}.


\medskip

The rest of the paper is organized as follows. In Section 2 we
discuss the basic volumetric properties of $2$-convex bodies. In
Section 3 we consider natural operations which preserve
$2$-convexity and its dual notion of $2$-smoothness, and prove
generalized versions of Proposition \ref{prop:Junge-And-Schatten}.
Section 4 treats various positions of $2$-convex bodies and their
interrelations. Section 5 deals with Gaussian marginals.
Throughout the text, we denote by $c, C, c^{\prime}$ etc. some
positive universal constant, whose value may change from line to
line. We will write $A \approx B$ to signify that $C_1 A \leq B
\leq C_2 A$ with universal constants $C_1,C_2>0$. We denote by
$D_n$ and $S^{n-1}$ the Euclidean unit ball and sphere in
$\Real^n$, respectively.

\medskip

\noindent \textbf{Acknowledgments.} Emanuel Milman would like to
sincerely thank his supervisor Prof. Gideon Schechtman for many
informative discussions.


\section{Volumetric properties} \label{sec:2} \label{sec:start}

Let $K \subset \RR^n$ be a centrally-symmetric convex body. Denote
by $\norm{\cdot}_K$ the norm whose unit ball is $K$. An
equivalent well-known characterization for $K$ to be $2$-convex
with constant $\alpha$ (e.g. \cite[Lemma 1.e.10]{LT-Book-II}) is
that for all $x,y \in \RR^n$:
\begin{equation} \label{eq:2-convex-equivalent}
\norm{x}_K^2 + \norm{y}_K^2 - 2 \norm{\frac{x+y}{2}}_K^2 \geq
\frac{\alpha'}{2} \norm{x-y}_K^2,
\end{equation}
where the relation between $\alpha$ and $\alpha'$ is summarized in
the following:
\begin{lem} \label{lem:equivalent-2-convex}
If $K$ is 2-convex with constant $\alpha$ then
(\ref{eq:2-convex-equivalent}) holds with $\alpha' = \alpha$. If
(\ref{eq:2-convex-equivalent}) holds for all $x,y \in \RR^n$,
then $K$ is 2-convex with constant $\alpha = \alpha' / 8$.
\end{lem}
It is also known (\cite{Nordlander}) that the Euclidean ball has
the best possible modulus of convexity, implying in particular
that $\alpha \leq 1/8$.

A basic observation due to Gromov and Milman
(\cite{Gromov-Milman}, see also \cite{ABV} for a simple proof) is
that if $K$ is uniformly convex with modulus of convexity
$\delta_K$, and $T \subset K$ with $|T| \geq \frac{1}{2} |K|$,
then for any $\eps
> 0$:
\begin{equation} \label{gm}
\frac{\left| (T + \eps K) \cap K \right|}{|K|} \geq 1 - 2 e^{-2 n
\delta_K(\eps)}.
\end{equation}
We will exploit (\ref{gm}) and obtain several interesting
consequences regarding mass distribution in $2$-convex sets. At
the heart of our argument is the following lemma, which is a
direct consequence of (\ref{gm}). We prefer to give a
self-contained proof, as this is a good opportunity to recreate
the elegant argument from \cite{ABV}. This lemma was also proved
in \cite{SchmuckyUniformlyConvexBodies}.

\begin{lemma} \label{lem:Psi-2-decay}
Let $K \subset \RR^n$ be a centrally-symmetric convex body. Assume
that $K$ is $2$-convex with constant $\alpha$, and that  $|K| =
1$. Fix $\theta \in S^{n-1}$ and denote $w = \sup_{x \in K}
|\langle x, \theta \rangle|$. Then for any $t > 0$:
\[
\text{Vol} \set{ x \in K ; \langle x, \theta \rangle > t } \leq 2
\exp \brac{-2 \alpha n (t / w)^2 }.
\]
\end{lemma}

\begin{proof}
Let $A(t) = \{ x \in K ; \langle x, \theta \rangle
> t \}$ and put $B = \{ x \in K ; \langle x, \theta \rangle < 0
\}$. Note that if $x \in A(t), y \in B$ then $\| x - y \|_K \geq
\frac{t}{w}$. According to the definition of $2$-convexity,
$$ \frac{B + A(t)}{2} \subset \left(1 - \alpha \left( \frac{t}{w} \right)^2 \right) K. $$
By the Brunn-Minkowski inequality,
\[
\sqrt{ |B| \cdot | A(t) |} \leq \abs{\frac{B + A(t)}{2}} \leq
\left( 1 - \alpha (t/w)^2 \right)^n \leq \exp{(- \alpha n
(t/w)^2)}.
\]
Since $|B| = 1/2$, we have:
\[
|A(t)| \leq 2 \exp{(- 2 \alpha n (t/w)^2)}.
\]
\end{proof}

Next, we present several consequences of Lemma
\ref{lem:Psi-2-decay}. The first one is the following observation.

\begin{lem} \label{lem1}
Let $K \subset \RR^n$ be a centrally-symmetric convex body. Assume
that $K$ is $2$-convex with constant $\alpha$ and volume 1, and
that $K$ is isotropic. Then:
\[
c \sqrt{\alpha} \sqrt{n} L_K D_n  \subset K,
\]
where $c > 0$ is a universal constant.
\end{lem}

\begin{proof}
Let $\theta \in S^{n-1}$ be arbitrary. For $t \in \Real$ set
\[
A(t) = K \cap \set{ x \in \RR^n ; \scalar{x, \theta} < t },
\]
and denote $f(t) = |A(t)|$. As before, we use $w = \sup_{x \in K}
|\langle x, \theta \rangle|$ to denote the width of $K$ in
direction $\theta$. By Lemma \ref{lem:Psi-2-decay}, we see that
for $t>0$:
\begin{equation} \label{eq:lem1-eq1}
f(t) \geq 1 - 2 \exp{(- 2 \alpha n (t/w)^2)}.
\end{equation}
On the other-hand, $f'(t) = \abs{K \cap \set{\scalar{x,\theta} =
t}}$ is a log-concave function by Brunn-Minkowski which is even,
and therefore attains its maximum at 0. Since $f'(0) \approx
1/L_K$ (e.g. \cite{Milman-Pajor-LK}), we see that:
\begin{equation} \label{eq:lem1-eq2}
f(t) \leq f(0) + t f'(0) \leq \frac{1}{2} + c \frac{t}{L_K}.
\end{equation}
Choosing $t = L_K / 4c$ and combining (\ref{eq:lem1-eq1}) and
(\ref{eq:lem1-eq2}), we see that $w \geq c' \sqrt{\alpha} \sqrt{n}
L_K$. Since the direction $\theta \in S^{n-1}$ was arbitrary, the
lemma follows.
\end{proof}

\medskip
Lemma \ref{lem1} entails Proposition \ref{prop:hyperplane} and
Proposition \ref{prop:finite-vr} at once. Indeed, since
$|\sqrt{n} D_n|^{1/n} \approx 1$ and $|K| = 1$, Lemma \ref{lem1}
implies that $L_K \leq c / \sqrt{\alpha}$. Proposition
\ref{prop:hyperplane} immediately follows (see, e.g.
\cite{Milman-Pajor-LK}). Since also $c < L_K$ (e.g.
\cite{Milman-Pajor-LK}), then Lemma \ref{lem1} implies that:
\[
c \sqrt{\alpha} \sqrt{n} D  \subset K,
\]
and Proposition \ref{prop:finite-vr} is established. Note that it
is quite unusual for a convex body to contain a large Euclidean
ball in isotropic position, even when the body has a small
volume-ratio. For instance, consider the convex body $K = \{ x \in
\RR^n; |x| \leq \sqrt{n}, |x_1| \leq 1 \}$, and let $\tilde{K}$ be
an isotropic linear image of $K$. It is easily seen that
$\tilde{K}$ does not contain a ball of radius larger than $c$,
although $K$ is isomrophic to a Euclidean ball, and clearly has a
finite volume-ratio.

\medskip

Another consequence of Lemma \ref{lem:Psi-2-decay} it the
following Proposition. As usual, the dual norm to
$\norm{\cdot}_K$ is defined by $\norm{x}^*_K = \sup_{y \in K}
\scalar{x,y}$, and its unit ball is called the polar body to $K$,
and denoted $K^\circ$. For $\theta \in S^{n-1}$, we define the
$\psi_2$-norm of the linear functional $\scalar{\cdot,\theta}$
w.r.t. the uniform measure on $K$ as:
\[
\norm{\scalar{\cdot, \theta}}_{L_{\psi_2(K)}} := \inf \set{
\lambda > 0 ; \frac{1}{|K|} \int_K e^{ \frac{\scalar{x,
\theta}^2}{\lambda^2}} dx \leq 2}.
\]
The $L_p$-norm is defined as:
\[
\norm{\scalar{\cdot, \theta}}_{L_p(K)} :=  \brac{\frac{1}{|K|}
\int_K \abs{\scalar{x, \theta}}^p dx}^{1/p}.
\]
It is well-known (e.g. \cite[Proposition 3.6]{JSZ}) that:
\[
\norm{\scalar{\cdot, \theta}}_{L_{\psi_2(K)}}
\approx \sup_{p\geq 2}
\frac{\norm{\scalar{\cdot,\theta}}_{L_p(K)}}{\sqrt{p}},
\]
implying in particular that:
\begin{equation} \label{eq:psi-2-and-diam}
\norm{\scalar{\cdot, \theta}}_{L_{\psi_2(K)}} \geq C
\frac{\norm{\theta}_K^*}{\sqrt{n}},
\end{equation}
since $\norm{\theta}_K^* \approx
\norm{\scalar{\cdot,\theta}}_{L_n(K)}$ (e.g.
\cite{Paouris-Small-Diameter}). Lemma \ref{lem:Psi-2-decay}
therefore implies:

\begin{prop} \label{prop:Psi-2-norm}
Let $K \subset \RR^n$ be a centrally-symmetric $2$-convex body
with constant $\alpha$. Then for all $\theta \in S^{n-1}$:
\[
C_1 \frac{\norm{\theta}_K^*}{\sqrt{n}} \leq \norm{\scalar{\cdot,
\theta}}_{L_{\psi_2(K)}} \leq C_2
\frac{\norm{\theta}_K^*}{\sqrt{\alpha} \sqrt{n}},
\]
where $C_1,C_2 > 0$ are two universal constants.
\end{prop}

Proposition \ref{prop:Psi-2-norm} provides us with a way to find
directions $\theta \in S^{n-1}$ for which $\VolSet{ x \in K ;
\scalar{x,\theta} \geq t }$ decays in a sub-gaussian rate, as
reflected by $\norm{\scalar{\cdot, \theta}}_{L_{\psi_2(K)}}$. As
a first application, note that for any convex body of volume one,
there exists a direction in which the width is smaller than $C
\sqrt{n}$ (otherwise the body would contain a Euclidean ball of
volume greater than one). Together with a straightforward
application of Markov's inequality, and denoting $M^*(K) =
\int_{S^{n-1}} \norm{\theta}^* d\sigma(\theta)$, we conclude the
following immediate corollary of Proposition
\ref{prop:Psi-2-norm}.

\begin{cor}
Let $K \subset \RR^n$ be a centrally-symmetric convex body. Assume
that $K$ is $2$-convex with constant $\alpha$ and volume 1. Then
there exists a universal constant $C>0$ such that:
\begin{enumerate}
\item
There exists a $\theta \in S^{n-1}$ such that:
\[
\norm{\scalar{\cdot, \theta}}_{L_{\psi_2(K)}} \leq C /
\sqrt{\alpha} \; .
\]
\item
\[
\sigma \set{ \theta \in S^{n-1} \; ; \;
\norm{\scalar{\cdot,\theta}}_{L_{\psi_2}(K)} \leq C
\frac{M^*(K)}{\sqrt{\alpha}\sqrt{n}} } \geq \frac{1}{2} \; .
\]
\end{enumerate}
\end{cor}

In Section 4, we will see several positions of a 2-convex body $K$
of volume 1 for which $M^*(K) \leq C \sqrt{n}$. The last
corollary implies that in these positions, at least half of the
directions have $\psi_2$-decay. We say that a body satisfying:
\[
\norm{\scalar{\cdot,\theta}}_{L_{\psi_2}(K)} \leq A \cdot
\vol{K}^{1/n}
\]
for all $\theta \in S^{n-1}$ is a $\psi_2$ body (with constant
$A$). In general, a $2$-convex body is not a $\psi_2$ body.
Indeed, as apparent from (\ref{eq:psi-2-and-diam}), a $\psi_2$
body (with constant $A$) of volume 1 always satisfies $diam(K)
\leq C A \sqrt{n}$, but any $l_p^n$ for $p<2$ (normalized to have
volume 1) already fails to satisfy this (with a universal
constant $A$) for large enough $n$. Here and henceforth,
$diam(K)$ denotes the diameter of $K$. Nevertheless, we can still
say the following:

\begin{prop}
Let $K \subset \RR^n$ be a centrally-symmetric convex body.
Assume that $K$ is $2$-convex with constant $\alpha$, has volume
1 and that it is isotropic. Then a random $\lfloor n/2 \rfloor$
dimensional section of $K$ is a $\psi_2$-body with high
probability.
\end{prop}
\begin{proof}
By definition, any section of $K$ is a $2$-convex body with the
same constant. By Proposition \ref{prop:finite-vr}, the isotropic
position is also a finite volume-ratio position for $K$, and $c
\sqrt{\alpha} \sqrt{n} D_n \subset K$. But by a classical result
of \cite{Szarek-fvr} and \cite{Szarek-TJ-fvr} (based on
\cite{Kashin-subspaces-of-l1}), a random $\lfloor n/2 \rfloor$
dimensional section $L \cap E$ of a convex body $L$ containing
$D_n$ is isomorphic to a Euclidean ball, and in particular
satisfies $diam(L \cap E) \leq C (|L|/|D_n|)^{2/n}$ with
probability greater than $1 - (1/2)^n$. Therefore:
\begin{equation} \label{eq:K-cap-E-inclusion}
c \sqrt{\alpha} \sqrt{n} (D_n \cap E) \subset K \cap E \subset
\frac{C'}{\sqrt{\alpha}} \sqrt{n} (D_n \cap E)
\end{equation}
with the same probability. Applying Proposition
\ref{prop:Psi-2-norm} to $K \cap E$ and using the left-hand-side
of (\ref{eq:K-cap-E-inclusion}) to compensate for the volume of $K
\cap E$, we see that:
\[
\norm{\scalar{\cdot, \theta}}_{L_{\psi_2(K \cap E)}} \leq
\frac{C'}{\alpha^{3/2}} \vol{K \cap E}^{2/n}
\]
for all $\theta \in S^{n-1} \cap E$. This concludes the proof.
\end{proof}


\section{Operations preserving 2-convexity}

We have already seen that (by definition) any section of a
$2$-convex body with constant $\alpha$ is itself a $2$-convex body
with the same constant. In this section we will consider several
additional natural operations which preserve $2$-convexity and
the dual notion of $2$-smoothness, and conclude with several new
results on the isotropic constant of different families of bodies.

The first natural operation to consider is taking projections.
Since this is the dual operation to taking sections, it will be
convenient to first introduce the dual notion to $2$-convexity,
which is $2$-smoothness. The modulus of smoothness of $K$ is
defined as the following function for $\tau>0$:
\begin{equation}
\rho_K(\tau) = \sup \set{ \frac{\norm{x+y}_K + \norm{x-y}_K}{2} -
1 \; ; \; \norm{x}_K \leq 1, \norm{y}_K \leq \tau }.
\end{equation}
A body $K$ is called ``$2$-smooth with constant $\beta$'' (see,
e.g. \cite[Chapter 1.e]{LT-Book-II}), if for all $\tau>0$:
\begin{equation}
\rho_K(\tau) \leq \beta \tau^2.
\end{equation}

It is well-known (e.g. \cite{LT-Book-II}) that the modulus of
smoothness is dual to the modulus of convexity (this can be
carefully formalized using Legendre transforms). We summarize
Propositions 1.e.2 and 1.e.6 from \cite{LT-Book-II} in the
following:
\begin{lem} \label{lem:duality}
Let $K$ be a centrally-symmetric convex body in $\Real^n$. Then
$K$ is $2$-convex with constant $\alpha$ iff $K^\circ$ is
$2$-smooth with constant $\frac{1}{16 \alpha}$.
\end{lem}

We will frequently refer to the Blaschke-Santalo inequality
(\cite{Santalo-inq}, the r.h.s. below) and its reverse form due to
Bourgain-Milman (\cite{Bourgain-Milman-vr-and-reverse-santalo},
the l.h.s. below), which together state that for any convex body
$K$:
\begin{equation} \label{eq:reverse-Santalo}
c \leq \nonumber \VolRatio{K}{D_n} \VolRatio{K^\circ}{D_n} \leq 1.
\end{equation}

Lemma \ref{lem:duality}, coupled with the Blaschke-Santalo
inequality or its reverse form, imply that we can translate many
volumetric results on $2$-convex bodies to $2$-smooth bodies. In
particular, Proposition \ref{prop:finite-vr} translates to the
fact that $2$-smooth bodies have finite \emph{outer-volume-ratio}.
We define the outer-volume-ratio of a body $K$ as:
\[
o.v.r.(K) = \inf_{\E \supset K} \left( \frac{|\E|}{|K|}
\right)^{\frac{1}{n}}
\]
where the infimum runs over all ellipsoids that contain $K$. If
$o.v.r.(K) < C$, for some universal constant $C>0$, it is
customary to say that $K$ has finite outer-volume-ratio. It is
well known (e.g. \cite{Milman-Pajor-LK}) that $L_K \leq C'
o.v.r.(K)$ for any convex body $K$. Combining everything
together, we have the following useful:

\begin{prop} \label{prop:2-smooth-ovr}
Let $K$ be a $2$-smooth convex body with constant $\beta$. Then
$o.v.r.(K) \leq C \sqrt{\beta}$. In particular, $L_K \leq C'
\sqrt{\beta}$.
\end{prop}

Note that if $K \subset T$ then $o.v.r.(K) \leq (|T|/|K|)^{1/n}
o.v.r.(T)$. The following is therefore an immediate corollary of
Proposition \ref{prop:2-smooth-ovr}:

\begin{cor} \label{cor:2-smooth-L_K}
Let $K$ be a centrally-symmetric convex body in $\Real^n$. Then:
\[
L_K \leq C \inf \set{ \left . \sqrt{\beta}
\brac{\frac{|T|}{|K|}}^{1/n} \; \right |
\begin{array}{c} K \subset T , \\ T \text{ is 2-smooth with constant $\beta$} \end{array} }
\]
\end{cor}

\medskip

We can now turn to investigate the action of taking projections
of $2$-convex and $2$-smooth bodies. For a subspace $E \subset
\RR^n$, we denote by $Proj_E$ the orthogonal projection onto $E$.
As evident from the definitions, any section of a $2$-smooth body
with constant $\beta$ is itself a $2$-smooth body with the same
constant. By passing to the polar body and using Lemma
\ref{lem:duality}, the duality between sections and projections
immediately implies:

\begin{lemma} \label{lem:projections}
Let $K \subset \RR^n$ be a $2$-convex ($2$-smooth) body with
constant $\gamma$. Then so is $Proj_E(K)$, with the same constant
$\gamma$, for any subspace $E \subset \RR^n$.
\end{lemma}

Using Lemma \ref{lem:projections}, a remarkable consequence of
Proposition \ref{prop:Psi-2-norm} is that the $\psi_2$-norm of the
linear functional $\scalar{\cdot,x}$ on a projection $Proj_E(K)$
of a $2$-convex body $K$, essentially depends (up to universal
constants) only on $x \in E$ and not on the subspace $E$. More
precisely:

\begin{prop}
Let $K \subset \RR^n$ be a $2$-convex body with constant
$\alpha$, and let $E$ be a $k$-dimensional subspace. Then for any
$x \in E$:
\[
C_1 \norm{x}_K^* \leq \norm{\scalar{\cdot,
x}}_{L_{\psi_2(Proj_E(K))}} \sqrt{k} \leq C_2
\frac{1}{\sqrt{\alpha}} \norm{x}_K^*
\]
\end{prop}
This is one of the rare cases where we can deduce volumetric
information on $Proj_E(K)$ from that of $K$. Typically, these two
bodies have different volumetric behaviour.

\smallskip

Let us consider other natural operations which preserve
2-convexity.
Unfortunately, the Minkowski sum is a bad candidate for this.
Indeed, even in $\RR^2$, the sum of two very narrow ellipsoids
which are perpendicular to each other, may be brought arbitrarily
close to a square, which is not 2-uniformly convex. Nevertheless,
there exists a well known natural summation operation, which
actually preserves both 2-uniform convexity and 2-uniform
smoothness. Recall that the \emph{2-Firey sum} of two convex
bodies $K$ and $T$, denoted by $K +_2 T$, is defined as the unit
ball of the norm satisfying:
\[
\norm{z}^2_{K+_2T} = \inf_{z = x+y} \norm{x}^2_K + \norm{y}^2_T.
\]
It is easy to see that the dual norms satisfy:
\[
(\norm{z}^*_{K+_2T})^2 = (\norm{z}^*_{K})^2 + (\norm{z}^*_{T})^2.
\]
We will refer to the latter operation as \emph{2-Firey
intersection}, and denote the 2-Firey intersection of $K$ and $T$
as $K\cap_2T$. Note that $(K\cap_2 T)^\circ = K^\circ +_2
T^\circ$.

\begin{lem} \label{lem:Firey}
Let $K$ and $T$ be 2-convex (smooth) bodies with constants
$\gamma_K$ and $\gamma_T$, respectively. Then so is their 2-Firey
sum $K+_2T$ and intersection $K \cap_2 T$, with constant
$\min\{\gamma_K,\gamma_T\} / 8$ ($\max\{\gamma_K,\gamma_T\} \cdot
8$).
\end{lem}
\begin{proof}
Obviously there is no loss in generality in assuming that
$\gamma_K = \gamma_T = \gamma$. Since $(K\cap_2 T)^\circ =
K^\circ +_2 T^\circ$, Lemma \ref{lem:duality} implies that the
case of 2-smooth bodies follows from the case of 2-convex bodies
by duality. We will therefore restrict ourselves to the latter
case, and assume that $K$ and $T$ are 2-convex with constant
$\gamma$. \\
By Lemma \ref{lem:equivalent-2-convex}, we have for $G = K,T$ and
for all $x,y\in \RR^n$:
\begin{equation} \label{eq:boaz2}
\norm{x}_G^2 + \norm{y}_G^2 - 2 \norm{\frac{x+y}{2}}_G^2 \geq
\frac{\gamma}{2} \norm{x-y}_G^2.
\end{equation}
Summing these two inequalities for $G=K$ and $G=T$, we see that
(\ref{eq:boaz2}) is also satisfied for $G = K\cap_2T$. Using Lemma
\ref{lem:equivalent-2-convex} again, this implies that $K\cap_2T$
is 2-convex with constant $\gamma/8$. \\
Next, for any $z_1,z_2 \in \RR^n$, write $z_i = x^K_i + x^T_i$ so
that:
\[
\norm{z_i}^2_{K+_2T} = \norm{x^K_i}^2_K + \norm{x^T_i}^2_T
\]
(by compactness the infimum is achieved). By Lemma
\ref{lem:equivalent-2-convex}, we know that for $G=K,T$:
\[
\norm{x^G_1}_G^2 + \norm{x^G_2}_G^2 \geq 2
\norm{\frac{x^G_1+x^G_2}{2}}_G^2 + \frac{\gamma}{2}
\norm{x^G_1-x^G_2}_G^2.
\]
Summing these two inequalities for $G=K$ and $G=T$ and denoting
$Z = K+_2T$, we have:
\begin{eqnarray}
\nonumber & & \norm{z_1}^2_Z + \norm{z_2}^2_Z = \norm{x^K_1}^2_K +
\norm{x^K_2}^2_K + \norm{x^T_1}^2_T + \norm{x^T_2}^2_T \\
\nonumber & \geq & 2 \brac{ \norm{\frac{x^K_1+x^K_2}{2}}_K^2 +
\norm{\frac{x^T_1+x^T_2}{2}}^2_T} + \frac{\gamma}{2} \brac{
\norm{x^K_1-x^K_2}_K^2 + \norm{x^T_1-x^T_2}_T^2} \\
\nonumber & \geq & 2 \norm{\frac{z_1+z_2}{2}}_Z^2 +
\frac{\gamma}{2} \norm{z_1-z_2}_Z^2,
\end{eqnarray}
where the last inequality follows from the definition of $Z =
K+_2T$ and the fact that $z_1 + z_2 = (x^K_1 + x^K_2) + (x^T_1 +
x^T_2)$ and $z_1 - z_2 = (x^K_1 - x^K_2) + (x^T_1 - x^T_2)$. Lemma
\ref{lem:equivalent-2-convex} implies that $K+_2T$ is 2-convex
with constant $\gamma / 8$.
\end{proof}

\begin{rem} \label{rem:Firey}
It is important to emphasize that the additional factor of 8
appearing in the Lemma is immaterial, and that the Lemma holds in
full generality when summing (intersecting) an arbitrary number of
bodies (with the same constant factor of 8).
\end{rem}

\medskip

We can now summarize our bounds for the isotropic constant in the
following statements. For a Banach space $X$, we denote by
$SQ_n(X)$ the class of unit balls of $n$-dimensional subspaces of
quotients of $X$. We denote $F^0_2SQ_n(X) = SQ_n(X)$, and by
induction:
\[
F^{i+1}_2SQ_n(X) = \set{ \bigwedge_{i=1}^l \bigoplus_{j=1}^{m_i}
K^i_j \; ; \; \set{K^i_j} \subset F^{i}_2SQ_n(X)},
\]
where $\bigwedge$ and $\bigoplus$ denote 2-Firey intersection and
sum, respectively. We set $F_2SQ_n(X) = \cup_{i=0}^\infty
F^i_2SQ_n(X)$. Note that it is possible to make the class
$F_2SQ_n(X)$ even richer, by alternately taking subspaces,
quotients, 2-Firey sums and 2-Firey intersections (since the
operation of 2-Firey sum is not distributive w.r.t. taking
subspace or 2-Firey intersection) starting from $X$, but this is a
complication which we wish to avoid. Lemmas \ref{lem:projections}
and \ref{lem:Firey}, together with Remark \ref{rem:Firey}, show
that if $X$ is 2-convex (2-smooth) with constant $\alpha$
($\beta$), then so is every member of $F_2SQ_n(X)$ with constant
$\alpha/8$ ($8 \beta$). Corollary \ref{cor:2-smooth-L_K}
therefore implies:

\begin{thm} \label{thm:2-smooth-L_K}
Let $K$ be a centrally-symmetric convex body in $\Real^n$, and
let $X$ be a 2-smooth Banach space with constant $\beta$. Then:
\[
L_K \leq C \sqrt{\beta} \inf \set{ \left .
\brac{\frac{|T|}{|K|}}^{1/n} \; \right | K \subset T , T \in
F_2SQ_n(X) }.
\]
\end{thm}

Consider $X = L_p$ for $2 \leq p < \infty$ in Theorem
\ref{thm:2-smooth-L_K}. Note that $X^* = L_q$ with $q = 1 +
1/(p-1)$, for which it is known (e.g. \cite[p. 63]{LT-Book-II})
that $X^*$ is 2-convex with constant equivalent to $1 / (p-1)$. By
Lemma \ref{lem:duality} this implies that $X$ is 2-smooth with
constant bounded by $C (p-1)$. We therefore have:

\begin{cor} \label{cor:L_p-L_K}
Let $K$ be a centrally-symmetric convex body in $\Real^n$. Then:
\[
L_K \leq C \inf \set{ \left . \sqrt{p}
\brac{\frac{|T|}{|K|}}^{1/n} \; \right | K \subset T , T \in
F_2SQ_n(L_p), p \geq 2 }.
\]
\end{cor}

This is a generalization of one half (the range $p \geq 2$) of a
Theorem of Junge
(\cite{Junge-slicing-problem-for-quotients-of-Lp}, see also
\cite{EMilman-DualMixedVolumes}):

\smallskip
\noindent\textbf{Theorem (Junge). }
\begin{equation} \label{eq:Junge}
\nonumber L_K \leq C \inf \set{ \left . \sqrt{p} \; q
\volRatio{T}{K} \: \right
| \; \begin{array}{c} K \subset T \; , \; T \in SQ_n(L_p) \; , \\
1<p<\infty \; , \; 1/p + 1/q = 1 \end{array} }.
\end{equation}
In fact, Junge showed that $L_p$ may be replaced by any Banach
space $X$ with finite type and bounded $gl_2(X)$ (the Gordon-Lewis
constant of $X$), in which case $\sqrt{p} \; q$ above should be
replaced by some constant depending on $X$.

We can also improve the second half of Junge's Theorem (in the
range $1<p\leq 2$) by replacing the factor of $q$ by $\sqrt{q}$.
Unfortunately, with our approach we have to insist that $K$ itself
is in $F_2SQ_n(L_p)$. Our version reads as follows:

\begin{thm} \label{thm:L_q-L_K}
Let $K \in F_2SQ_n(L_p)$ for $1<p\leq 2$, and let $q$ be given by
$1/p + 1/q = 1$. Then:
\[
L_K \leq C \sqrt{q}.
\]
\end{thm}

The latter is an immediate corollary of the the fact that $L_p$
for $1<p\leq 2$ is 2-convex with constant equivalent to $p-1$
(e.g. \cite[Chapter 1.e]{LT-Book-II}), combined with the following
general Theorem, which is a consequence of Proposition
\ref{prop:hyperplane}:

\begin{thm} \label{thm:2-convex-L_K}
Let $X$ be a 2-convex Banach space with constant $\alpha$, and let
$K \in F_2SQ_n(X)$. Then:
\[
L_K \leq C \frac{1}{\sqrt{\alpha}}.
\]
\end{thm}

Another interesting example is obtained by taking $X$ to be the
space of all $m$ by $m$ complex or real matrices, equipped with
the norm $\norm{A} = (tr (A A^*)^{p/2})^{1/p}$, the so-called
$l_p$-Schatten-Class which will be denoted by $S_p^m$.
It was observed in \cite{Schatten-LK-Bounded} that the isotropic
constants of these spaces are uniformly bounded (in $m$), which
is especially interesting in the range $1\leq p<2$, since for $p
\geq 2$ it is known that the unit ball of $S_p^m$ (or any of its
subspaces) has finite outer volume-ratio. In the former range, it
has been recently shown in \cite{Guedon-Paouris-Schatten} that (in
particular) the isotropic constants of several special subspaces
of $S_p^m$ are also uniformly bounded. Although our method does
not extend to $p=1$, we can show the following result, which in
particular demonstrates that the same is true for any subspace of
quotient of $S_p^m$, provided that $p$ is bounded away from 1.
The modulus of convexity (and smoothness) of $S_p^m$ was
estimated by N. Tomczak-Jaegermann in \cite{TJ-Schatten-Modulus},
where it was shown that $\delta_{S_p^m} \approx \delta_{L_p}$. It
follows that $S_p^m$ is 2-convex with constant equivalent to
$p-1$ for $1<p\leq 2$, which together with Theorem
\ref{thm:2-convex-L_K} gives:

\begin{thm} \label{thm:Schatten}
Let $K \in F_2SQ_n(S_p^m)$ for $1 < p \leq 2$ and $m\geq n$, and
let $q$ be given by $1/p + 1/q = 1$. Then:
\[
L_K \leq C \sqrt{q}.
\]
\end{thm}

It is clear that the case $p=1$ in Theorem \ref{thm:L_q-L_K} and
Theorem \ref{thm:Schatten} must serve as a break-down point for
our method. Indeed, since $S_1^m$ contains $l_1^m$ as a subspace
(of the diagonal matrices), and since every convex body may be
approximated as the unit ball of a quotient of $l_1^m$ for
large-enough $m$, or simply as the quotient of $L_1$, a similar
result for $p=1$ in either theorem would solve the Slicing
Problem.



\section{Equivalence between positions of $2$-convex bodies}

For the results of this section, we will need to recall a few
basic notions from Banach space theory. The (Rademacher) type-$p$
constant of a Banach space $X$ (for $1 \leq p \leq 2$), denoted
$T_p(X)$, is the minimal $T>0$ for which:
\[
\brac{\mathbb{E} \snorm{\sum_{i=1}^m \eps_i x_i}^2}^{1/2} \leq T
\brac{\sum_{i=1}^m \norm{x_i}^2}^{1/2}
\]
for any $m \geq 1$ and any $x_1,\ldots,x_m \in X$, where
$\set{\eps_i}$ are independent, identically distributed random
variables uniformly distributed on $\set{-1,1}$ and $\mathbb{E}$
denotes expectation. Similarly, the cotype-$q$ constant of $X$
(for $2 \leq q \leq \infty$), denoted $C_q(X)$, is the minimal
$C>0$ for which:
\[
\brac{ \mathbb{E} \snorm{\sum_{i=1}^m \eps_i x_i}^2}^{1/2} \geq
\frac{1}{C} \brac{\sum_{i=1}^m \norm{x_i}^q}^{1/q}
\]
for any $m \geq 1$ and $x_1,\ldots,x_m \in X$. We say that $X$
\emph{has} type $p$ (cotype $q$) if $T_p(X) < \infty$ ($C_q(X) <
\infty$). We also say that $X$ is \emph{of} type $p$ (cotype $q$)
if $p = \sup \set{ p' ; X \text{ has type } p' }$ ($q = \inf \set{
q' ; X \text{ has cotype } q' }$).

Let $L_2(\set{-1,1}^m,X)$ denote the space of $X$-valued
functions on the discrete cube $\set{-1,1}^m$, equipped with the
norm $(\mathbb{E} \norm{f(\eps_1,\ldots,\eps_m)}^2 )^{1/2}$. We
denote by $Rad_m(X)$ the Rademacher projection on
$L_2(\set{-1,1}^m,X)$ (see \cite{Milman-Schechtman-Book}), and
denote $\norm{Rad(X)} = sup_m \norm{Rad_m(X)}$ where
$\norm{Rad_m(X)}$ is the operator norm of $Rad_m(X)$. By duality,
it is easy to verify that $\norm{Rad(X^*)} = \norm{Rad(X)}$, and
it is clear that $\norm{Rad_m(X)} = sup_{E \subset X}
\norm{Rad_m(E)}$ where the supremum runs over all
finite-dimensional subspaces of $X$.

One of the most important results in the so-called local-theory
of Banach spaces is a theorem by Pisier who showed that
$\norm{Rad(X)}$ may be bounded from above by an (explicit)
function of $T_p(X)$ when $p>1$, concluding that $\norm{Rad(X)} <
\infty$ when $X$ has type $p>1$. When $p=2$, there is a much
easier argument, going back to a remark at the end of the work by
Maurey and Pisier \cite{Maurey-Pisier-Theorem} (see also
\cite[Remark 2.11]{Georgian-K-Convexity-Paper} for an explicit
proof), showing (without any constants!):

\begin{lem} \label{lem:Rad-type-2}
$\norm{Rad(X)} \leq T_2(X).$
\end{lem}

The next lemma, which gives a non-quantitive estimate of the
opposite inequality (for the general $p$ case) using a
compactness argument, is a known consequence of the Maurey-Pisier
Theorem \cite{Maurey-Pisier-Theorem}:

\begin{lem} \label{lem:Maurey-Pisier}
There exists a function $C(R) : \Real_+ \rightarrow \Real_+$ such
that any finite-dimensional Banach space $X$ with $\norm{Rad(X)}
\leq R$ satisfies $T_{p(R)}(X) \leq C(R)$ with $p(R) =  1 + 1 /
C(R)$.
\end{lem}
\begin{proof}[Sketch of proof]
Assume that this is not true for some $R>0$. This means that there
exist finite-dimensional Banach spaces $X_i$ with $\norm{Rad(X_i)}
\leq R$ and $T_{1+1/i}(X_i) > i$. The latter easily implies that
$dim(X_i) \rightarrow \infty$, since always $T_p(X_i) \leq
T_2(X_i) \leq \sqrt{dim(X_i)}$ for any $1 \leq p \leq 2$ ($X_i$
is $\sqrt{dim(X_i)}$-isomorphic to a Hilbert space $H_i$ by John's
Theorem, and $T_2(H_i)=1$). We now construct an infinite
dimensional Banach space $X$ as the $l_2$ sum of the $X_i$'s,
i.e. for $x = (x_i)_{i \geq 1}$ with $x_i \in X_i$ define
$\norm{x}_X = (\sum_{i \geq 1} \norm{x_i}^2_{X_i})^{\frac{1}{2}}$
and set $X = \set{ x ; \norm{x}_X < \infty }$ endowed with the
norm $\norm{\cdot}_X$. It is elementary to check that
$\norm{Rad(X)} \leq R$, and since $X$ contains each $X_i$ as a
subspace we must have that $X$ is of type 1. The latter implies
by the Maurey-Pisier Theorem (actually we only need the type 1
case, which is due to Pisier \cite{Pisier-Type-1-Spaces}) that
$X$ contains $(1+\epsilon)$ isometric copies of $l_1^m$ for
arbitrary $\epsilon >0$ and $m$, and as a consequence
$\norm{Rad(X)} \geq sup_m \norm{Rad(l_1^m)} = \infty$. We arrive
to a contradiction, so the assertion is proved.
\end{proof}

Let us return to the study of $2$-convex bodies. We recall the
following classical result (e.g. \cite[Theorem
1.e.16]{LT-Book-II}). For completeness, we sketch the proof.
\begin{lem} \label{lem:type-cotype-2}
\hfill
\begin{enumerate}
\item
Let $K$ be a $2$-convex body with constant $\alpha$. Then
$C_2(X_K) \leq \frac{C}{\sqrt{\alpha}}$.
\item
Let $K$ be a $2$-smooth body with constant $\beta$. Then $T_2(X_K)
\leq C \sqrt{\beta}$.
\end{enumerate}
\end{lem}
\begin{proof}
(1) easily follows from the equivalent characterization
(\ref{eq:2-convex-equivalent}) of a $2$-convex body, which
asserts that for any $x_1,x_2 \in \Real^n$:
\[
\mathbb{E} \norm{\eps_1 x_1 + x_2}^2 = \frac{1}{2}(\norm{x_2 +
x_1}^2 + \norm{x_2 - x_1}^2) \geq \alpha \norm{x_1}^2 +
\norm{x_2}^2.
\]
Hence by induction, since $\alpha < 1$:
\[
\mathbb{E} \norm{\sum_{i=1}^m \eps_i x_i}^2 \geq \alpha
\sum_{i=1}^m \norm{x_i}^2,
\]
for any $x_1,\ldots,x_m \in \Real^m$, which concludes the proof
of (1) (even without a constant!). (2) follows either by duality
or similarly from the equivalent characterization of a $2$-smooth
body (e.g. \cite[Theorem A.7]{BL-Book}):
\[
\norm{x+y}^2 + \norm{x-y}^2 - 2\norm{x}^2 \leq C\beta \norm{y}^2,
\]
for every $x,y \in \Real^n$.
\end{proof}

We are now ready to conclude the following useful:

\begin{lem} \label{lem:main-Rad}
Let $K$ be a $2$-convex body with constant $\alpha$. Then:
\begin{enumerate}
\item
\[
\norm{Rad(X_K)} \leq C / \sqrt{\alpha}.
\]
\item
There exists a $p>1$ which depends on $\alpha$ only, such that:
\[
T_p(X_K) \leq 1 / (p-1).
\]
\end{enumerate}
\end{lem}
\begin{proof}
By Lemma \ref{lem:duality}, $K^\circ$ is 2-smooth with constant
$1/(16\alpha)$, and so by Lemmas \ref{lem:Rad-type-2} and
\ref{lem:type-cotype-2} we see that:
\[
\norm{Rad(X)} = \norm{Rad(X^*)} \leq T_2(X^*) \leq
\frac{C}{\sqrt{\alpha}},
\]
which concludes the proof of (1). Applying Lemma
\ref{lem:Maurey-Pisier}, we immediately deduce (2).
\end{proof}

Lemmas \ref{lem:type-cotype-2} and \ref{lem:main-Rad} allow us to
deduce several interesting results about $2$-convex bodies. By a
classical result of Figiel and Tomzcak-Jaegermann on the
$l$-position (\cite{l-position}), for any convex body $K$ there
exists a position for which $M(K) M^*(K) \leq C \norm{Rad(X_K)}$,
and in fact this is satisfied in the \emph{minimal mean-width
position}. The latter is defined (up to orthogonal rotations) as
the volume-preserving affine image of $K$ for which $M^*(K)$ is
minimal. Recall that we always have:
\begin{equation} \label{eq:Urysohn}
\frac{1}{M(K)} \leq \VolRad{K} \leq M^*(K),
\end{equation}
where $\VolRad{K} = (\vol{K} / \vol{D_n})^{1/n}$ and the first
inequality follows from Jensen's inequality while the second is
Urysohn's inequality. We therefore deduce that in the minimal
mean-width position, a $2$-convex body $K$ with constant $\alpha$
satisfies:
\begin{equation} \label{eq:boaz}
M^*(K) \leq \frac{C}{\sqrt{\alpha}} \VolRad{K},
\end{equation}
which is essentially the best possible by (\ref{eq:Urysohn}). We
will refer to (\ref{eq:boaz}) as ``$M^*(K)$ is bounded'', omitting
the reference to the volume-radius. As we shall see, there are
many advantages of working with a position in which $M^*(K)$ is
bounded.

Our next Proposition shows that whenever we have a good upper
bound on $M^*(K)$, $K$ is essentially isotropic. For convenience,
we define $M_2^*(K) = (\int_{S^{n-1}} (\norm{\theta}^*_K)^2
d\sigma(\theta))^{1/2}$, which is well known to be equivalent to
$M^*(K)$ (by Kahane's inequality for instance).

\begin{proposition} \label{prop:essential-isotropic}
For any $2$-convex body $K$ with constant $\alpha$ and volume 1,
we have:
\[
\int_K \abs{x} dx \leq C \frac{M^*(K)}{\sqrt{\alpha}}.
\]
\end{proposition}
\begin{proof}

\begin{eqnarray}
\nonumber \int_K \abs{x} dx \leq \brac{\int_K \abs{x}^2 dx}^{1/2}
= \sqrt{n} \brac{\int_K \int_{S^{n-1}} \scalar{x,\theta}^2
d\sigma(\theta) dx}^{1/2}  \\
\nonumber = \sqrt{n} \brac{\int_{S^{n-1}} \int_K
\scalar{x,\theta}^2 dx \; d\sigma(\theta)}^{1/2} = \sqrt{n}
\brac{\int_{S^{n-1}} \norm{\scalar{\cdot,\theta}}^2_{L_2(K)}
d\sigma(\theta)}^{1/2} \\
\nonumber \leq C \sqrt{n} \brac{\int_{S^{n-1}}
\norm{\scalar{\cdot,\theta}}^2_{L_{\psi_2}(K)}
d\sigma(\theta)}^{1/2} \leq \frac{C'}{\sqrt{\alpha}}
\brac{\int_{S^{n-1}} (\norm{\theta}_K^*)^2 d\sigma(\theta)}^{1/2},
\end{eqnarray}
where we used Proposition \ref{prop:Psi-2-norm} in the last
inequality. The last term is equal to $\frac{C'}{\sqrt{\alpha}}
M_2^*(K)$, which is majorized by $\frac{C''}{\sqrt{\alpha}}
M^*(K)$.
\end{proof}

The last Proposition has an interesting consequence regarding
2-Firey sums of 2-convex bodies in minimal mean-width position, or
in any bounded $M^*$ position in general.
\begin{corollary}
Let $K$ and $T$ be 2-uniformly convex bodies, such that $M^*_2(K)
\leq C_K \VolRad{K}$ and $M^*_2(T) \leq C_T \VolRad{T}$ (and
therefore essentially isotropic). Then $M^*_2(K+_2T) \leq
\max(C_K,C_T) \VolRad{K+_2T}$. In particular, $K+_2T$ is
essentially isotropic.
\end{corollary}
\begin{proof}
Notice that $(M^*_2)^2$ is clearly additive with respect to
2-Firey sums, whereas by \cite{Lutwak-Firey-Sums}
$\vol{K+_2T}^{2/n} \geq \vol{K}^{2/n} + \vol{T}^{2/n}$. The claim
then easily follows.
\end{proof}

\medskip

An additional property of any position for which $M^*(K)$ is
bounded, is that it automatically satisfies half of the
conditions of being in a 2-regular M-position. Recall that a
convex body $K$ in $\RR^n$ is said to be in $a$-regular
M-position ($0<a\leq 2$) if its homothetic copy $K'$, normalized
to that $\vol{K'} = \vol{D_n}$, satisfies:
\begin{equation} \label{eq:regular-position}
N(K',t D_n) \leq \exp(C n / t^a) \text{ and } N((K')^\circ,t D_n)
\leq \exp(C n / t^a),
\end{equation}
for $t \geq 1$, where $N(K,L)$ is the covering number of $K$ by
$L$ (see \cite{GiannopoulosMilmanHandbook}) and $C>0$ is a
universal constant. It was shown by Pisier
(\cite{Pisier-Regular-M-Position}) that an $a$-regular M-position
for $0<a<2$ always exists (with a constant $C$ in
(\ref{eq:regular-position}) depending only on $a$). When $M^*(K)$
is bounded and $\vol{K} = \vol{D_n}$, by Sudakov's inequality
(\cite{GiannopoulosMilmanHandbook}):
\[
N(K,tD_n) \leq \exp(C n (M^*(K) / t)^2) \leq \exp(C n / t^2)
\]
for $t\geq 1$, so half of the condition for being in a
$2$-regular M-position is satisfied. In general, the other half of
the condition, namely:
\begin{equation} \label{eq:dual-M-position}
N(K^\circ,tD_n) \leq \exp(C n / t^2),
\end{equation}
does not follow from knowing that $M^*(K)$ is bounded.
Nevertheless, we mention two cases where this would follow. If
$K$ is in minimal mean-width position and $\vol{K} = \vol{D_n}$,
in which case both $M(K)$ and $M^*(K)$ are bounded by
(\ref{eq:Urysohn}), then (\ref{eq:dual-M-position}) follows from
Sudakov's inequality applied to $K^\circ$. Another case is when
$K$ is in a finite volume-ratio position with bounded $M^*(K)$
(remember that we know that $K$ has finite volume-ratio), in
which case (\ref{eq:dual-M-position}) is trivially satisfied. The
second case, if it exists, will be preferred over the first,
since it adds the finite-volume ratio position property (which is
not guaranteed in general by the minimal mean-width position), in
particular implying that $M(K)$ is bounded.

Luckily, for a $2$-convex body, there exists an "all-in-one"
position which gives all of the above mentioned properties:
bounded $M^*$, having finite volume-ratio (and therefore being in
a 2-regular M-position) and essential isotropicity. This position
is exactly John's maximal-volume ellipsoid position. This follows
from the following useful lemma from
\cite{EMilman-DualMixedVolumes} (which appeared first in an
equivalent form in \cite{Davis-etal-Lemma}):
\begin{lemma} \label{lem:type-2-lemma}
For any convex body $K$ in John's maximal-volume ellipsoid
position, the following holds:
\[
M_2^*(K) b(K) \leq T_2(X_K^*),
\]
where $b(K) = \max_{\theta \in S^{n-1}} \norm{\theta}_K$.
\end{lemma}
\noindent For a 2-convex body $K$ with constant $\alpha$, the
polar body is 2-smooth with constant $1/(16 \alpha)$, and
therefore by Lemma \ref{lem:type-cotype-2}, $X_K^*$ has type 2
with constant $T_2(X_K^*) \leq C/\sqrt{\alpha}$. Noting that
$M^*(K) \leq M_2^*(K)$, Lemma \ref{lem:type-2-lemma} therefore
gives:
\begin{corollary}
A 2-convex body $K$ with constant $\alpha$ in John's
maximal-volume ellipsoid position, satisfies:
\[
M^*(K) b(K) \leq \frac{C}{\sqrt{\alpha}}.
\]
\end{corollary}

Since $M^*(K) b(K)$ is invariant under homothety, we may assume
above that $\vol{K} = \vol{D_n}$, in which case $b(K) \geq 1$ (by
volume consideration) and $M^*(K) \geq 1$ (by Urysohn's
inequality). We therefore see that in John's maximal-volume
ellipsoid position $M^*(K) \leq C/\sqrt{\alpha} \VolRad{K}$. The
similar bound on $b$ implies again that $K$ has finite-volume
ratio, $v.r.(K) \leq C/\sqrt{\alpha}$, with the same bound (up to
a possible constant) as in Proposition \ref{prop:finite-vr}.
Proposition \ref{prop:essential-isotropic} coupled with the latter
bound on $M^*(K)$ in John's position, imply Proposition
\ref{prop:inverse-finite-vr} stated in the Introduction.

\medskip

One last additional property that we would like our "all-in-one"
position to satisfy is having a small-diameter: if $\vol{K} =
\vol{D_n}$, we would like to have $diam(K) \leq C (n /
\log{n})^{1/2}$. The motivation for this requirement comes from
\cite{ABP}, where it was shown that if an isotropic 2-convex body
has small-diameter in the above sense, then most of its marginals
are approximately Gaussian (see \cite{ABP} or Section
\ref{sec:Gaussian} for more details). It is easy to check that
this requirement is indeed satisfied by all the $l_p^n$ unit
balls for $1<p\leq 2$ (normalized to have the appropriate volume).

Unfortunately, the small-diameter requirement is not satisfied
for a general 2-convex body
in isotropic position, as illustrated by the
following:

\begin{example} \label{example}
Let:
\[
T = \left \{ (x,y) \in \RR^2 ; x^2 + \left( |y| + 1 \right)^2
\leq 2 \right \}.
\]
The set $T$ is $2$-convex with constant $c$, and has two ``cusps",
at $(1,0)$ and $(-1,0)$. Denote by $K \subset \RR^n$ the
revolution body of $T$ around the $y$-axis, namely:
\[
K = \left \{ (x_1,...,x_n) \in \Real^n; \left( (x_1^2 + \ldots +
x_{n-1}^2)^{1/2}, x_n \right) \in T \right \}.
\]
It is easy to check that $K$ is $2$-convex with constant $c'$. Let
$\widetilde{K} \subset \RR^n$ be an isotropic image of $K$ of
volume 1. Then $diam(\widetilde{K}) \geq c'' n$.
\end{example}
\begin{proof}[Sketch of proof]
Around its "cusp" hyperplane $e_n^\perp$, $K$ looks like a
two-sided cone, and therefore half of the volume of $K$ lies
inside the slab $\{ x\in\Real^n ; \abs{\scalar{x,e_n} } \leq
c(n)/n \}$ with $c(n) \approx 1$. But in isotropic position of
volume 1, half of the volume of $\widetilde{K}$ lies inside slabs
of width in the order of $L_K$ (and $L_K \approx 1$ by
Proposition \ref{prop:hyperplane}). This means that we must
inflate $K$ by an order of $n$ in the direction of $e_n$ when
passing to $\widetilde{K}$, implying that $diam(\widetilde{K})
\geq c'' n$.
\end{proof}

Nevertheless, the following proposition shows that in L\"{o}wner's
\emph{minimal}-volume outer ellipsoid position, the small-diameter
requirement is satisfied, although we are not able to guarantee
any of the other "good" properties satisfied by John's
maximal-volume ellipsoid position. We note that $K$ is in
L\"{o}wner's position iff $K^\circ$ is in John's position.

\begin{proposition} \label{prop:small-diam}
Let $K$ be any 2-convex body with constant $\alpha$ and volume 1.
Then there exists a constant $\lambda>0$ which depends on $\alpha$
only, such that in L\"{o}wner's minimal-volume outer ellipsoid
position, $diam(K) \leq \frac{C}{\lambda} n^{1/2-\lambda}$.
\end{proposition}
\begin{proof}
Apply Lemma \ref{lem:type-2-lemma} to $K^\circ$, which by duality
is in John's maximal-volume ellipsoid position. Then:
\[
M_2(K) diam(K) \leq T_2(X_K).
\]
Since $M_2(K) \geq \VolRad{K}^{-1} = 1$ by Jensen's inequality, it
is enough to show that $T_2(X_K)$ is bounded by $C
n^{1/2-\lambda}$. By Lemma \ref{lem:main-Rad}, we know that there
exists a $p>1$ which depends on $\alpha$ only, such that
$T_p(X_K) \leq 1/(p-1)$, so it remains to pass from type-$p$ to
type-2. But this is an easy consequence of a result by
Tomczak-Jaegermann (\cite{TJ-Type-Cotype-fev-vectors}), who
showed that it is enough to evaluate the type 2 constant of an
$n$-dimensional Banach space on $n$ vectors. If $x_1 \ldots x_n$
is any sequence in $\RR^n$, then by H\"{o}lder's inequality:
\[
\mathbb{E} \snorm{\sum_{i=1}^n \eps_i x_i}_K \leq \frac{1}{p-1}
\brac{\sum_{i=1}^n \norm{x_i}^p_K }^ {\frac{1}{p}} \leq
\frac{n^{\frac{1}{p} - \frac{1}{2}}}{p-1} \brac{\sum_{i=1}^n
\norm{x_i}^2_K }^ {\frac{1}{2}}.
\]
Therefore $T_2(X_K) \leq \frac{C}{\lambda} n^{1/2 - \lambda}$,
for $\lambda = 1 - 1/p$.
\end{proof}

We conclude this section by mentioning that the results of
Section \ref{sec:2} imply that for $2$-convex bodies, the
isotropic position is a 1-regular M-position. Indeed, since the
isotropic position is also a finite volume-ratio position, the
second half of condition (\ref{eq:regular-position}) is trivially
satisfied. The first half is satisfied by the result from
(\cite{Hartzoulaki-PhD} or \cite[Proposition
5.4]{Klartag-Isomophic-Slicing}), which shows that this is always
the case for any isotropic body for which $L_K$ is bounded. Note
that \cite[Theorem 5.6]{GiannopoulosMilmanMeanWidth} (which uses
Dudley's entropy bound) enables us to bound the mean-width of a
convex body in an $a$-regular M-position, which for a $1$-regular
position gives:
\[
M^*(K) \leq C diam(K)^{1/2} \VolRad{K}^{1/2}.
\]
Since $diam(K) \leq C \sqrt{n} L_K \VolRad{K}$ in isotropic
position (e.g. \cite{Milman-Pajor-LK}), we conclude that $M^*(K)
\leq C(\alpha) n^{1/4} \VolRad{K}$ for any $2$-convex body $K$
with constant $\alpha$ in isotropic position. It is still unclear
to us whether the isotropic position is always a 2-regular
M-position, which would imply (as above) that $M^*(K) \leq
C(\alpha) \log(n) \VolRad{K}$.

\medskip

To summarize, we have seen the following implications for a
2-convex body:
\begin{itemize}
\item
Minimal mean-width position implies essential isotropicity and a
2-regular M-position.
\item
John's maximal-volume ellipsoid position implies finite
volume-ratio position, essential minimal mean-width, 2-regular
M-position and essential isotropicity.
\item
L\"{o}wner's minimal-volume outer ellipsoid position implies
"small-diameter".
\item
Isotropic position implies finite volume-ratio position and
1-regular M-position.
\end{itemize}


\section{Gaussian marginals} \label{sec:Gaussian}

Similarly to the $2$-convex case, we say that a convex body $K$
is $p$-convex (with constant $\alpha$) if its modulus of convexity
satisfies $\delta_K(\epsilon) \geq \alpha \epsilon^p$ for all
$\epsilon \in (0,2)$. Let us also denote $d_K = diam(K)$. It is
well-known and easy to see (e.g. \cite{Ledoux-Book} or follow the
argument in Lemma \ref{lem:Psi-2-decay}) that the Gromov-Milman
Theorem (\ref{gm}) immediately implies the following:
\begin{lem} \label{lem:p-convex-concentration}
Let $K$ be a $p$-convex body with constant $\alpha$ and of volume
1. For any 1-Lipschitz function $f$ on $K$ denote by $Med(f)$ the
median of $f$, i.e. the value for which $\VolSet{ x \in K ; f(x)
\geq Med(f) } \geq 1/2$ and $\VolSet{ x \in K ; f(x) \leq Med(f) }
\geq 1/2$. Then:
\[
\VolSet{ x \in K ; f(x) \geq Med(f) + t} \leq 2 \exp(-2 \alpha n
(t / d_K)^{p}).
\]
\end{lem}
Let us denote $E(f) = \int_K f(x) dx$. As in \cite{ABP}, we deduce
from Lemma \ref{lem:p-convex-concentration} that $\abs{E(f) -
Med(f)} \leq C \; d_K (\alpha n) ^{-\frac{1}{p}}$. We therefore
have:
\[
\VolSet{ x \in K ; \abs{ f(x) - E(f) } \geq t + C d_K (\alpha
n)^{-\frac{1}{p}} } \leq 4 \exp(-2 \alpha n
\brac{\frac{t}{d_K}}^{p}),
\]
and it is easy to check that this implies:
\begin{lem} \label{lem:p-convex-concentration2}
With the same notations as in Lemma
\ref{lem:p-convex-concentration}:
\[
\VolSet{ x \in K ; \abs{ f(x) - E(f) } \geq t } \leq 4 \exp(- 2
c^p \alpha n \brac{\frac{t}{d_K}}^{p}).
\]
\end{lem}
Using this, it was shown in \cite{ABP} that if $K$ is an isotropic
$p$-convex body (with constant $\alpha$) with $\vol{K} = 1$ and
$diam(K) \leq R \sqrt{n}$, then:
\[
\VolSet{ x\in K ; \abs{ \frac{|x|}{\sqrt{n}} - L_K } \geq R t }
\leq 4 \exp( - 2 c^p \alpha n t^p ).
\]
Choosing $t = C (\frac{\log(n)}{\alpha n})^{1/p}$, this implies:
\begin{equation} \label{eq:p-convex-concentration2}
\VolSet{ x\in K ; \abs{ \frac{|x|}{\sqrt{n}} - L_K } \geq C R
\brac{\frac{\log(n)}{\alpha n}}^{1/p} } \leq \frac{1}{n}.
\end{equation}
The authors of \cite{ABP} conclude that if $R \ll (\alpha n /
\log(n))^{1/p}$, (\ref{eq:p-convex-concentration2}) implies a
concentration of the volume of $K$ inside a spherical shell
around a radius of $\sqrt{n} L_K$. It was shown in \cite{ABP} that
such a concentration implies that most marginals of the uniform
distribution on $K$ will have an approximately Gaussian
distribution (see Theorem \ref{thm:extended-ABP} below).
Unfortunately, our investigation of the case $p=2$ shows that
this condition on $R$ is not satisfied in general by isotropic
$2$-convex bodies, as demonstrated by Example \ref{example}.
Nevertheless, Proposition \ref{prop:small-diam} shows that in
L\"{o}wner's minimal-volume ellipsoid position, we do have $R
\leq C n^{1/2 - \lambda}/\lambda$ where $\lambda$ depends only on
the $2$-convexity constant of $K$. In this case, the concentration
result of \cite{ABP} still holds, with the minor change that $L_K$
in (\ref{eq:p-convex-concentration2}) is replaced by $\int_K
\abs{x} dx / \sqrt{n}$ (note that this value is always greater
than $c_1 L_K \geq c_2$, e.g. \cite{Milman-Pajor-LK}). Although
$K$ is no longer isotropic, it is possible to generalize the
argument in \cite{ABP} to a body in arbitrary position. This is
done in \cite{EMilman-Gaussian-Marginals}, where the following is
shown:

\begin{thm}[Generalization of \cite{ABP}] \label{thm:extended-ABP}
Let $K$ be a centrally-symmetric convex body in $\Real^n$ of
volume 1, and assume that for some $\rho>0$ and $\epsilon < 1/2$:
\begin{equation} \label{eq:extended-ABP-condition}
\VolSet{ x \in K ; \abs{\frac{|x|}{\sqrt{n}} - \rho} \geq \epsilon
\rho } \leq \epsilon.
\end{equation}
For $\theta \in S^{n-1}$ denote $g_\theta(s) = \Vol{K \cap \set{s
\theta + \theta^\perp}}$ and let $\rho_\theta^2 =
\int_{-\infty}^\infty s^2 g_\theta(s) ds$. Denote the Gaussian
density with variance $\rho^2$ by $\phi(s) = \frac{1}{\sqrt{2 \pi}
\rho} \exp(-\frac{s^2}{2 \rho^2})$ and let $H(\theta) = \sup_{t>0}
\abs{\int_{-t}^t g_\theta(s) ds - \int_{-t}^t \phi(s) ds}$. Then
for any $0<\delta<c$:
\begin{eqnarray}
\nonumber & & \!\!\!\!\!\!\!\!\!\!\!\!\!\!\!\!\! \sigma \set{
\theta \in S^{n-1} ; H(\theta) \leq \delta + 4 \epsilon + \frac{c_1}{\sqrt{n}} } \\
\label{eq:ABP-gaussians} &\geq & 1- C_1 C_{iso}(K) \sqrt{n}
\log{n} \exp\brac{- \frac{c_2 n \delta^2}{C_{iso}(K)^2}},
\end{eqnarray}
where:
\[
\rho_{max} = \max_{\theta \in S^{n-1}} \rho_\theta \; , \;
\rho_{avg} = \int_{S^{n-1}} \rho_\theta d\sigma(\theta) \; , \;
C_{iso}(K) = \frac{\rho_{max}}{\rho_{avg}}.
\]
\end{thm}
\begin{rem}
As usual, it is easy to verify that $\rho_{avg}$ and $\rho$ above
are equivalent to within absolute constants (since $\epsilon <
1/2$).
\end{rem}

If $T$ is a volume preserving linear transformation such that
$\widetilde{K} = T(K)$ is isotropic, then clearly $\rho_{max} =
\norm{T^{-1}}_{op} L_K$, where $\norm{\cdot}_{op}$ denotes the
operator norm. Since $\rho^2_{avg} \approx \frac{1}{n} \int_K
|x|^2 dx \geq L_K^2$ (e.g. \cite{Milman-Pajor-LK}), it follows
that $C_{iso}(K) \leq C \norm{T^{-1}}_{op}$. Hence, knowing that
$ r D_n \subset \widetilde{K}$ and $K \subset R D_n$ would imply
that $C_{iso}(K) \leq C R / r$. By Lemma \ref{lem1} and
Proposition \ref{prop:small-diam}, $c \sqrt{\alpha} \sqrt{n} L_K
D_n \subset \widetilde{K}$ and $K \subset C n^{1-\lambda} /
\lambda$ in L\"{o}wner's position, where $\lambda>0$ depends only
on $\alpha$. We therefore have in this position:
\[
C_{iso}(K) \leq \min(\frac{ C n^{1/2 - \lambda} }{\sqrt{\alpha}
\lambda L_K}, C \sqrt{n}).
\]
Hence, regardless of its a-priori diameter, by putting a 2-convex
body $K$ with constant $\alpha$ in L\"{o}wner's position, we
deduce by Proposition \ref{prop:small-diam}, Lemma
\ref{lem:p-convex-concentration2} and Theorem
\ref{thm:extended-ABP} that most marginals of $K$ are
approximately Gaussian in the above sense, where the level of
proximity ($\epsilon$ above) depends only on $\alpha$.
Summarizing, we have:

\begin{thm} \label{thm:2-convex-Gaussian}
Let $K$ be a 2-convex body with constant $\alpha$ and volume 1.
Assume that $K$ is in L\"{o}wner's minimal-volume outer ellipsoid
position. Then with the same notations as in Theorem
\ref{thm:extended-ABP} and with $\rho = \int_K |x| dx /
\sqrt{n}$, we have for any $0 < \delta < c$:
\[
\sigma\set{ \theta \in S^{n-1} ; H(\theta) \leq \delta + 4
\epsilon + \frac{c_1}{\sqrt{n}}} \geq 1-  n^{5/2} \exp\brac{-c_2
\alpha n^{2\lambda} \lambda^2 \delta^2},
\]
where $\epsilon = C \sqrt{\log n} \: \alpha^{-1/2} \lambda^{-1}
n^{-\lambda}$ and $\lambda = \lambda(\alpha) > 0$ depends on
$\alpha$ only.
\end{thm}

Before concluding, we remark that placing a 2-convex body $K$ in
L\"{o}wner's position is just a convenient "pre-processing" step.
In fact, in any position we always have at least one approximately
Gaussian marginal (in the above sense); it just happens that in
L\"{o}wner's position we can show this for "most" marginals w.r.t.
the Haar probability measure on the unit sphere, and this would
equally be true in an arbitrary position by choosing a different
measure (the one induced by the change of positions, for
example). The reason is that the metric given by $H(\theta)$ in
Theorem \ref{thm:extended-ABP} is invariant under
volume-preserving linear transformations. More precisely, given
such a $T$, and any body $K$ and $\rho>0$, it is immediate to
check that:
\[
\int_{-t}^t (g^K_\theta(s) - \phi_{\rho}(s)) ds =
\int_{-\frac{t}{|T(\theta)|}}^{\frac{t}{|T(\theta)|}} (
g^{T(K)}_{\frac{T(\theta)}{|T(\theta)|}}(s) - \phi_{\rho
|T(\theta)| }(s)) ds,
\]
so by Theorem \ref{thm:2-convex-Gaussian} we can control the
supremum over $t>0$ of either expressions for at least one
$\theta \in S^{n-1}$ if $K$ is a 2-convex body in L\"{o}wner's
position and $\rho = \int_K |x| dx / \sqrt{n}$.

\bibliographystyle{amsalpha}
\bibliography{../../ConvexBib}

\end{document}